\documentclass[12pt]{article}

\usepackage{graphicx}
\usepackage{epstopdf}

\usepackage[T1]{fontenc}
\usepackage[utf8]{inputenc}
\usepackage{authblk}

\usepackage[numbers]{natbib}
\usepackage[OT4]{fontenc}
\usepackage[utf8]{inputenc}

\newcommand{\Exp}{{\rm I\hspace{-0.8mm}E}}

\newcommand{\Var}{{\bf Var}}

\newcommand{\iz}{{\rm \rlap Z\kern 2.2pt Z}}

\newcommand{\proof}{\noindent {\bf Proof:} \ }

\newtheorem{theorem}{Theorem}

\newtheorem{proposition}{Proposition}

\newtheorem{remark}{Remark}

\author[1]{Zbigniew Michna\footnote{Corresponding author\\ Email: zbigniew.michna@ue.wroc.pl\\ Tel/fax: +48713680335}}
\author[2]{Peter Nielsen}
\affil[1]{Department of Mathematics and Cybernetics

Wroc{\l}aw University of Economics}
\affil[2]{Department of Mechanical and Manufacturing Engineering Aalborg University}

\title{\bf\LARGE {The impact of lead time forecasting on the bullwhip effect}}

\date{}

\begin{document}

\maketitle

\bibliographystyle{abbrv}

\begin{abstract}
In this article we quantify the bullwhip effect (the variance amplification in replenishment orders)  when demands and lead times are predicted in a simple two-stage supply chain with one supplier and one retailer. In recent research the impact of stochastic order lead time on the bullwhip effect is investigated, but the effect of needing to predict / estimate the lead time is not considered in the supply chain models.
Under uncertainty conditions it is necessary to estimate the lead time for a member of the supply chain to place an order.
We find a new cause of the bullwhip effect in the form of lead time forecasting and we give an exact form of the bullwhip effect measure (the ratio of variances) when demands and lead times are predicted by moving averages. In the bullwhip effect measure we discover two terms amplifying the effect which are the result of lead time estimation. 

\vspace{5mm}
\end{abstract}

\section{Introduction}
The bullwhip effect was recognized by Forrester \cite{fo:58} in the middle of the twentieth century and was coined as a term by Procter \& Gamble management.
This phenomenon appears in supply chains as the variance amplification in replenishment orders if one moves up in
a supply chain (see Disney and Towill \cite{di:to:03} and Geary et al. \cite{ge:di:to:06} for the definition and historical review). It is considered harmful because of its consequences which are (see e.g. Buchmeister et al. \cite{bu:pa:pa:po:08}):
excessive inventory investment, poor customer service level, lost revenue, reduced productivity, more difficult decision-making, sub-optimal transportation, sub-optimal production etc. This makes it critical to find the root causes of the bullwhip effect and to quantify the increase in demands variability at each stage of the supply chain as this is directly linked to costs. In the current state of research typically five main causes of the bullwhip effect are considered (see e.g. Lee et al. \cite{le:pa:wh:97a} and \cite{le:pa:wh:97b}): demand forecasting, non-zero lead time, supply shortage, order batching and price fluctuation. To decrease the variance amplification in a supply chain (i.e. to reduce the bullwhip effect) we need to identify all factors causing the bullwhip effect and to quantify their impact on the effect.

In this research we will investigate another cause of the bullwhip effect that is {\bf lead time forecasting}. It is well known from inventory theory that the mean and variability of lead time of a supplier affects the
inventory and order decisions of its customer. Although lead times are typically considered deterministic, they are actually not in many supply chains (see Chatfield et al. \cite{ch:ki:ha:ha:04}).
The impact of stochastic lead times in inventory
systems has been intensively studied in the literature see e.g. Bagchi et al. \cite{ba:ha:ch:86},
Hariharan and Zipkin \cite{ha:zi:95}, Mohebbi and
Posner \cite{mo:po:98}, Song \cite{so:94a} and \cite{so:94b}, Song and Zipkin
\cite{so:zi:93} and \cite{so:zi:96} and Zipkin \cite{zi:86}.
Recent research investigates the influence of stochastic lead times on the bullwhip effect (see e.g. So and Zheng \cite{so:zh:03}, Duc et al. \cite{du:lu:ki:08} or Kim et al. \cite{ki:ch:ha:ha:06} and references therein) but none of the current research addresses the consequences of inherent need to estimate / forecast lead time when it is stochastic. These works investigate the impact of random stochastic lead times on the bullwhip effect through the characteristics of their distribution e.g. mean value or variance. In the paper So and Zheng \cite{so:zh:03} supplier's delivery lead time depends on the existing order backlog at
the supplier which means that it is not deterministic but it depends on the retailer's order quantities.
They solve this problem numerically.

In the typical approach if one assumes that a certain feature is random there is a need to predict its value for the next periods.
In our situation the relationship between supplier's lead time and its customer (a retailer) order quantities is very strong, especially when the supplier operates at tight capacity and has difficulties to adjust
capacity and maintain constant delivery lead
time to its customers. We should also notice that the retailer
order quantities can in turn determine the
delivery time performance of the supplier.
If a retailer observes uncertainty in demands and lead times (i.e. they are random) and he 
wants to place an order to a supplier due to a certain stock policy to fulfill customer orders in a
timely manner, he needs to predict future customer's demands and future supplier's lead times.  In other words the retailer needs to project his costumer's future demands over 
his supplier's lead time to determine the
appropriate order quantity to this supplier. This is done by the so-called lead time demand forecasting
to have the necessary required inventory to meet customer demands over the lead time. Lead time demand forecasting
can be executed by demand forecasting and lead time forecasting.
Thus we can not avoid that the value of a future lead time is necessary to determine an order quantity to the supplier. It yields a need to predict lead times based on their previous values. Practically a retailer needs to estimate 
(to forecast) the value of the next lead time to make an order to a supplier. This need for lead time estimation has not been noticed
in previous works on the impact of a stochastic lead time on the bullwhip effect.

In this paper we find a new cause
of the bullwhip effect which is lead time estimation through forecasting and we quantify its impact on the variance amplification in replenishment orders. Many papers assuming a deterministic lead time have studied the influence of different methods of demand forecasting on the bullwhip effect such as simple moving average, exponential
smoothing, and minimum-mean-squared-error forecasts when demands are independent identically distributed or constitute integrated moving-average, autoregressive process or autoregressive-moving average (see Graves \cite{gr:99}, Lee et al. \cite{le:so:ta:00}, Chen et al. \cite{ch:dr:ry:si:00a} and \cite{ch:dr:ry:si:00b}, Alwan et al. \cite {al:li:ya:03}, Zhang \cite{zh:04} and Duc et al. \cite{du:lu:ki:08}). Also the recent works on the impact of lead times on the bullwhip effect by Agrawal et al. \cite{ag:se:sh:09} and Li and Liu \cite{li:li:13} should be noted. However, of these papers the first one does not consider stochastic lead times and the second one investigates a transition state model with uncertainties in demands, production process, supply chain structure, inventory policy implementation and especially vendor order placement lead time delays. They find
a maximally allowable vendor order placement lead time delay such that the supply chain system is exponentially stabilizable.
This approach uses dynamical control systems theory and is not probabilistic (for similar models see the references in Li and Liu \cite{li:li:13}).

In this paper we consider moving averages as methods of demand and lead time forecasting and we find an exact form of the bullwhip effect measure related to the prediction of lead times and demands. 
More precisely we investigate a model where:
\begin{enumerate}
\item[a)] a supply chain contains two stages and consists of a retailer who receives client demands and a supplier (customers $\leftrightarrow$ retailer  $\leftrightarrow$ supplier (manufacturer));
\item[b)] customer demands constitute an iid sequence;
\item[c)] lead times between the supplier and the retailer constitute an iid sequence;
\item[d)] the retailer uses the order-up-to level policy to make an order to the supplier;
\item[e)] the retailer predicts the future values of demands and the future value of lead times
based on the simple moving average method using past observations that is we propose the following lead time demand forecast
$$
\widehat{D_t^L}=\sum_{i=0}^{\widehat{L_t}-1}\widehat{D}_{t+i}\,,
$$
where $\widehat{L_t}$ is the forecast for a next lead time of the order made at the beginning of a period $t$ and
$\widehat{D}_{t+i}$ denotes the forecast for a demand for the period $t+i$ at the beginning of a period $t$.
\end{enumerate}
The crucial point of our approach is the last subpoint e) and differs from the previous approaches. Namely 
in the work of Duc et al. \cite{du:lu:ki:08}
the lead time demand forecast is defined as follows
$$
\widehat{D_t^L}=\sum_{i=0}^{L_t-1}\widehat{D}_{t+i}\,,
$$
where $L_t$ is the next lead time at the beginning of a time $t$.  The value of $L_t$ the retailer does not know at the beginning of the time $t$ when he makes an order to the supplier. This means that the last lead time demand forecasting is not feasible in practice. 
The paper of Kim et al. \cite{ki:ch:ha:ha:06} also investigates a stochastic lead time in supply chains 
and proposes lead time demand forecasting. More precisely the simply moving average method for lead time demand is proposed that is
\begin{equation}\label{kim1}
\widehat{D_t^L}=\frac{1}{p}\sum_{j=1}^{p}D_{t-j}^L\,,
\end{equation}
where $p$ is the delay parameter of the prediction and $D_{t-j}^L$ is the previous known lead time demand of the order made at the beginning of the time $t-j$. This approach is practically feasible. Let us notice that 
\begin{equation}\label{kim2}
D_{t-j}^L=\sum_{i=0}^{L_{t-j}-1}D_{t-j+i}\,,
\end{equation}
where $L_{t-j}$ is a lead time of an order made at the beginning of the time $t-j$ and $D_{t-j+i}$ is the demand from the period $t-j+i$. Combining 
(\ref{kim1}) and (\ref{kim2}) we get a double sum and we can not exchange the sums because $L_{t-j}$ are different (compare it with Kim et al. \cite{ki:ch:ha:ha:06}).

In our approach we show that the bullwhip effect measure contains two summands depending on 
lead time forecasting. These terms amplify the value of the bullwhip effect measure and are the evidence that
lead time estimation in itself  is another cause of the bullwhip effect. 

\section{Supply chain model}
We will model a supply chain with two stages that is one retailer and one supplier. In our approach to the problem of lead time forecasting we assume that the retailer observes demands $D_t$ of his costumers (usually $t$ denotes a time period and $D_t$ is a demand during a period of the same length). More precisely we will assume that $\{D_t\}_{t=-\infty}^\infty$ constitutes a sequence of independent identically distributed random variables with $\Exp D_t=\mu_D$ and $\Var D_t=\sigma_D^2$ and a generic random variable for demands will be denoted by $D$. Similarly lead times are introduced that is $L_t$ is the lead time for an order placed by the retailer to the supplier at the beginning of the period $t$. Random variables of lead times $\{L_t\}_{t=-\infty}^\infty$ are independent identically distributed with $\Exp L_t=\mu_L$ and $\Var L_t=\sigma_L^2$ and a generic random variable for lead times we will denote by $L$.  Let us note that we do not impose any assumptions on the distributions of $D$ and $L$. We assume only that their second moments are finite. The sequences $\{D_t\}_{t=-\infty}^\infty$ and $\{L_t\}_{t=-\infty}^\infty$ are independent of each other. The lead time demand at the beginning of a period $t$ is defined as follows
\begin{equation}\label{ltd}
D_t^L=D_t+D_{t+1}+\dots...+D_{t+L_t-1}=
\sum_{i=0}^{L_t-1}D_{t+i}\,.
\end{equation}
This value is not known for the retailer at the beginning of a period $t$ but he needs to forecast its value to make an order to the supplier. The natural way to do this is to predict demands and lead times. 
If $\widehat{D}_{t+i}$ denotes the forecast for a demand for the period $t+i$ at the beginning of a period $t$ (that is after $i+1$ periods, $i=0,1,\ldots$) and ${\cal{F^D}}_{t-1}=\sigma(D_{t-1}, D_{t-2},\ldots)$ is the sigma algebra generated by the demands up to a time $t-1$ then 
$$
\widehat{D}_{t+i}\in {\cal{F^D}}_{t-1}\,, 
$$
which means that the forecast at the beginning of the period $t$ for a period $t+i$ is a function of the previous known demands $\{D_{t-1}, D_{t-2},\ldots\}$. Similarly since lead times are random the retailer needs to predict their values for the next periods to make an order. Let ${\cal{F^L}}_{t-1}=\sigma(L_{t-1}, L_{t-2},\ldots)$ be the sigma algebra generated by lead times up to a time $t-1$. Thus if  $\widehat{L_t}$
is the forecast for a next lead time at the beginning of a period $t$ then generally 
$$
\widehat{L_t}\in {\cal{F^L}}_{t-1}\,, 
$$
which means that the forecast for a period $t$ is a function of the previous known lead times $\{L_{t-1}, L_{t-2},\ldots\}$. Thus the retailer making an order to a supplier puts the following forecast for a lead time demand as follows
\begin{equation}\label{eltd0}
\widehat{D_t^L}=\sum_{i=0}^{\widehat{L_t}-1}\widehat{D}_{t+i}\,.
\end{equation}
Employing
the moving average forecast method with the length $n\geq 1$ for demand forecasting we get
\begin{equation}\label{df}
\widehat{D}_{t+j}=\frac{1}{n}\sum_{i=1}^n D_{t-i}\,, 
\end{equation}
where $j=0,1,\ldots$ and $D_{t-i}$ $i=1,2,\ldots, n$ are demands which have been observed by the retailer till the
beginning of a period $t$. Similarly, the retailer predicts a lead time. Precisely, using
the moving average forecast method with the length $m\geq 1$ for lead time forecasting we obtain
\begin{equation}\label{ltf}
\widehat{L_t}=\frac{1}{m}\sum_{i=1}^m L_{t-i}\,, 
\end{equation}
where $L_{t-i}$ $i=1,2,\ldots, m$ are lead times which have been observed by the retailer till the
beginning of a period $t$. If we want to be more precise we need to assume that the distribution
of lead is such that 
$$
L_t\leq M
$$
where $M>0$ that is lead times are bounded by $M$. This we assume to avoid the situation that
for example the lead time $L_{t-1}$ is not known at the beginning of the time $t$ when we make an order. Then lead time forecasting is the following
\begin{equation}\label{ltfm}
\widehat{L_t}=\frac{1}{m}\sum_{i=1}^m L_{t-M-i}
\end{equation}
that is we get back at least $M$ periods. For simplicity we will use in our calculation the lead time forecast given in (\ref{ltf}) because one can see slightly modifying the proof of Th. \ref{bmmt} that the bullwhip effect measure is the same under
assumption that lead times are bounded and applying the lead time forecast given in (\ref{ltfm}).
Thus by eq. (\ref{eltd0}), (\ref{df}) and (\ref{ltf}) we get the forecast for a lead time demand as follows
\begin{equation}\label{eltd}
\widehat{D_t^L}=\widehat{L_t}\widehat{D_t}=
\frac{1}{mn}\sum_{i=1}^m L_{t-i}\sum_{i=1}^n D_{t-i}\,.
\end{equation}
We have to indicate that a similar idea to (\ref{eltd}) appeared in Chatfield at el. \cite{ch:ki:ha:ha:04}
but there the bullwhip effect measure is simulated without showing the realtion between lead time forecasting and the bullwhip effect.
We can employ the lead time forecast (\ref{ltfm}) to (\ref{eltd}) but as we mentioned this does not
affect the bullwhip effect measure.
Moreover in our model the retailer applies a base stock policy that is a simple order-up-to level inventory policy. Let $S_t$ be the inventory position at the beginning of a period $t$ (later an order is placed).
If the order-up-to level policy is employed then $S_t$ is determined in the following way
\begin{equation}\label{st}
S_t=\widehat{D_t^L}+z\widehat{\sigma_t}\,,
\end{equation}
where 
$$
\widehat{\sigma_t}^2=\Var(D_t^L-\widehat{D_t^L})
$$
is the variance of the forecast error for the lead time demand
and $z$ is the normal z-score that specifies the probability that demand is fulfilled by the on-hand inventory and it can be found based on a given service level. 
In some articles $\widehat{\sigma_t}^2$ is defined more practically that is instead of variance it is taken the empirical variance of  $D_t^L-\widehat{D_t^L}$. This complicates calculations very much but we must mention that the estimation of  $\widehat{\sigma_t}^2$ increases the size of the bullwhip effect. These two approaches coincide if z=0.
Thus the order quantity $q_t$ placed
at the beginning of a period $t$ is
\begin{equation}\label{qt}
q_t=S_t-S_{t-1}+D_{t-1}\,.
\end{equation}
Our main purpose is to find $\Var q_t$ and then to calculate the following bullwhip effect measure
$$
BM=\frac{\Var q_t}{\Var D_t}\,.
$$
\begin{proposition}
The variance of the forecast error for the lead time demand does not depend on $t$ and is as follows
\begin{eqnarray*}
\lefteqn{\widehat{\sigma_t}^2=\Var(D_t^L-\widehat{D_t^L})}\\
&&=\mu_L\sigma_D^2+\frac{\sigma_L^2\mu_D^2(m+1)}{m}+\frac{\mu_L^2\sigma_D^2}{n}
+\frac{\sigma_L^2\sigma_D^2}{mn}\,.
\end{eqnarray*}
\end{proposition}
\proof
By the eq. (\ref{ltd}) and (\ref{eltd}) and assuming independence, we get that $\Exp D_t^L=\Exp \widehat{D_t^L}= \mu_L\mu_D$ and
\begin{eqnarray*}
		\lefteqn{\Var(D_t^L-\widehat{D_t^L})=\Exp(D_t^L-\widehat{D_t^L})^2}\\
&&=\Exp\left(\sum_{i=0}^{L_t-1}D_{t+i}\right)^2+\frac{1}{m^2n^2}\Exp\left(\sum_{i=1}^m L_{t-i}\right)^2\Exp\left(\sum_{i=1}^n D_{t-i}\right)^2\\
&&\,\,\,\,\,\,\,\,-2\Exp\left(\sum_{i=0}^{L_t-1}D_{t+i}\right)\frac{1}{m}\Exp\left(\sum_{i=1}^m L_{t-i}\right)
\frac{1}{n}\Exp\left(\sum_{i=1}^n D_{t-i}\right)\\
&&=\Exp L\Exp D^2 +\Exp(L(L-1))(\Exp D)^2\\
&&\,\,\,\,+\frac{1}{m^2n^2}[m\Exp L^2+m(m-1)(\Exp L)^2]
[n\Exp D^2+n(n-1)(\Exp D)^2]\\
&&\,\,\,\,-2(\Exp L\Exp D)^2\\
&& = \mu_L(\sigma_D^2+\mu_D^2)+(\sigma_L^2+\mu_L^2-\mu_L)\mu_D^2\\
&&\,\,\,\,+\frac{1}{mn}(\sigma_L^2+m\mu_L^2)(\sigma_D^2+n\mu_D^2)-2\mu_L^2\mu_D^2\\
&&=\mu_L\sigma_D^2+\frac{\sigma_L^2\mu_D^2(m+1)}{m}+\frac{\mu_L^2\sigma_D^2}{n}
+\frac{\sigma_L^2\sigma_D^2}{mn}
\end{eqnarray*}
which finishes the proof.

Since the variance of the forecast error for the lead time demand is independent of $t$ we have from the eq. (\ref{st}) and (\ref{qt})
$$
q_t=\widehat{D_t^L}-\widehat{D_{t-1}^L}+D_{t-1}
$$
which permits to calculate the variance of $q_t$.
\begin{proposition}
The variance of an order quantity in a period $t$ is given as
$$
\Var q_t=\frac{2\sigma_L^2\sigma_D^2(m+n-1)}{m^2n^2}+\frac{2\sigma_L^2\mu_D^2}{m^2}+\frac{2\mu_L^2\sigma_D^2}{n^2}+\frac{2\mu_L\sigma_D^2}{n}+\sigma_D^2\,.
$$
\end{proposition}
\proof
Let us note that
\begin{eqnarray*}
\widehat{D_{t-1}^L}&=&\frac{1}{mn}\sum_{i=1}^m L_{t-1-i}\sum_{i=1}^n D_{t-1-i}\\
&=&\frac{1}{mn}\sum_{i=2}^{m+1} L_{t-i}\sum_{i=2}^{n+1} D_{t-i}\\
&=&\frac{1}{mn}\left(\sum_{i=1}^{m}L_{t-i}+L_{t-m-1}-L_{t-1}\right)\left(\sum_{i=1}^{n} D_{t-i}+D_{t-n-1}-D_{t-1}\right)\\
&=&\widehat{D_t^L}+\frac{1}{mn}(D_{t-n-1}-D_{t-1})\sum_{i=1}^m L_{t-i}+
\frac{1}{mn}(L_{t-m-1}-L_{t-1})\sum_{i=1}^n D_{t-i}\\
&&\,\,\,\,+\frac{1}{mn}(L_{t-m-1}-L_{t-1})(D_{t-n-1}-D_{t-1})\,.
\end{eqnarray*}
Thus we get
\begin{eqnarray*}
q_t&=&-\frac{1}{mn}(D_{t-n-1}-D_{t-1})\sum_{i=1}^m L_{t-i}-
\frac{1}{mn}(L_{t-m-1}-L_{t-1})\sum_{i=1}^n D_{t-i}\\
&&\,\,\,\,\,-\frac{1}{mn}(L_{t-m-1}-L_{t-1})(D_{t-n-1}-D_{t-1})+D_{t-1}\,.
\end{eqnarray*}
By independence it is easy to notice that $\Exp q_t=\mu_D$. So let us compute the second moment
of $q_t$
\begin{eqnarray*}
\lefteqn{\Exp q_t^2=}\\
&&=\frac{1}{m^2n^2}\Exp(D_{t-n-1}-D_{t-1})^2\Exp(\sum_{i=1}^m L_{t-i})^2\\
&&\,\,\,\,+\frac{1}{m^2n^2}
\Exp(L_{t-m-1}-L_{t-1})^2\Exp(\sum_{i=1}^n D_{t-i})^2\\
&&\,\,\,\,+\frac{1}{m^2n^2}\Exp(L_{t-m-1}-L_{t-1})^2\Exp(D_{t-n-1}-D_{t-1})^2+
\Exp D_{t-1}^2\\
&&\,\,\,\,+\frac{2}{m^2n^2}\Exp[(L_{t-m-1}-L_{t-1})\sum_{i=1}^m L_{t-i}]
\Exp[(D_{t-n-1}-D_{t-1})\sum_{i=1}^n D_{t-i}]\\
&&\,\,\,\,+\frac{2}{m^2n^2}\Exp(D_{t-n-1}-D_{t-1})^2\Exp[(L_{t-m-1}-L_{t-1})\sum_{i=1}^m L_{t-i}]\\
&&\,\,\,\,-\frac{2}{mn}\Exp[D_{t-1}(D_{t-n-1}-D_{t-1})]\Exp(\sum_{i=1}^m L_{t-i})\\
&&\,\,\,\,+\frac{2}{m^2n^2}\Exp[(L_{t-m-1}-L_{t-1})^2\Exp[(D_{t-n-1}-D_{t-1})\sum_{i=1}^n D_{t-i}]\\
&&\,\,\,\,-\frac{2}{mn}\Exp(L_{t-m-1}-L_{t-1})\Exp(D_{t-1}\sum_{i=1}^n D_{t-i})\\
&&\,\,\,\,-\frac{2}{mn}\Exp(L_{t-m-1}-L_{t-1})\Exp[D_{t-1}(D_{t-n-1}-D_{t-1}) ]\\
&&=\frac{2\sigma_D^2}{m^2n^2}[m(\sigma_L^2+\mu_L^2)+m(m-1)\mu_L^2]+\frac{2\sigma_L^2}{m^2n^2}[n(\sigma_D^2+\mu_D^2)+n(n-1)\mu_D^2]\\
&&\,\,\,\,+\frac{4\sigma_L^2\sigma_D^2}{m^2n^2}+\sigma_D^2+\mu_D^2
+\frac{2}{m^2n^2}(\mu_L^2-\sigma_L^2-\mu_L^2)(\mu_D^2-\sigma_D^2-\mu_D^2)\\
&&\,\,\,\,+\frac{4\sigma_D^2}{m^2n^2}(\mu_L^2-\sigma_L^2-\mu_L^2)
-\frac{2}{mn}(\mu_D^2-\sigma_D^2-\mu_D^2)m\mu_L\\
&&\,\,\,\,+\frac{4\sigma_L^2}{m^2n^2}(\mu_D^2-\sigma_D^2-\mu_D^2)-0-0\\
&&=\frac{2\sigma_L^2\sigma_D^2(m+n-1)}{m^2n^2}+\frac{2\sigma_L^2\mu_D^2}{m^2}+\frac{2\mu_L^2\sigma_D^2}{n^2}+\frac{2\mu_L\sigma_D^2}{n}+\sigma_D^2+\mu_D^2
\end{eqnarray*}
which gives the thesis.

Thus we can derive the exact form of the bullwhip effect measure.
\begin{theorem}\label{bmmt}
The measure of the bullwhip effect has the following form
$$
BM=\frac{\Var q_t}{\Var D_t}=
\frac{2\sigma_L^2(m+n-1)}{m^2n^2}+\frac{2\sigma_L^2\mu_D^2}{m^2\sigma_D^2}+\frac{2\mu_L^2}{n^2}+\frac{2\mu_L}{n}+1\,.
$$
\end{theorem}
\begin{remark}
We get the same formula if we employ the lead time forecast (\ref{ltfm}) under assumption that
lead times are bounded.
\end{remark}
Let us analyze the formula. The first summand in the formula includes the impact of the forecast of
lead times and demands. The second summand shows the influence of the prediction of lead times.
The third and fourth ones give the amplification of the variance by demand forecasting.
The effect is very large (see the next section and the tables below) if we take $m=1$ that is in the case if the forecast of a next lead time is based on one last observation of the lead time then we get
\begin{eqnarray*}
BM=\frac{\Var q_t}{\Var D_t}&=&
\frac{2\sigma_L^2}{n}+\frac{2\sigma_L^2\mu_D^2}{\sigma_D^2}+\frac{2\mu_L^2}{n^2}+\frac{2\mu_L}{n}+1\\
&=&\frac{2\sigma_L^2\mu_D^2}{\sigma_D^2}+\frac{2\mu_L^2}{n^2}+\frac{2(\mu_L+\sigma_L^2)}{n}+1\,.
\end{eqnarray*}
If lead times are deterministic that is $L_t=L=const.$ then the bullwhip effect is described by
$$
BM=\frac{\Var q_t}{\Var D_t}=
\frac{2L^2}{n^2}+\frac{2L}{n}+1\,,
$$
which is consistent with the result of Chen et al. \cite{ch:dr:ry:si:00a}. We should notice that Duc et al.
\cite{du:lu:ki:08} also obtained the result of Chen et al. \cite{ch:dr:ry:si:00a} in a special case and
as an exact value of the bullwhip effect (not a lower bound). Chen et al. \cite{ch:dr:ry:si:00a} get this as a lower bound because they define the error $\widehat{\sigma_t}$ as the empirical variance of  $D_t^L-\widehat{D_t^L}$.

Now we investigate what happens if the number of past observations of lead times or demands are large
that is if $m\rightarrow\infty$ or $n\rightarrow\infty$. So if the number of past lead times included in the forecast (the delay parameter of forecasting) goes to infinity we get
$$
\lim_{m\rightarrow\infty}BM=
\frac{2\mu_L^2}{n^2}+\frac{2\mu_L}{n}+1\,.
$$
This shows that the impact of the prediction of lead times disappears if the number of previous lead
times included in the forecast is very large. 
Similarly if the number of demands 
used in the prediction is growing to infinity then
$$
\lim_{n\rightarrow\infty}BM=\frac{2\sigma_L^2\mu_D^2}{m^2\sigma_D^2}+1\,.
$$
The effect has not disappeared and it remains constant if the ratio $\mu_D^2/\sigma_D^2$ does not
change and it is linear with respect to $\sigma_L^2$. Moreover this term can be very harmful if $m$ is small (see the next section and the tables below).

\section{Numerical examples}
We will numerically investigate the measure of the bullwhip effect. Especially we will consider every term in formula
of Th. \ref{bmmt}. Thus let us put
\begin{equation}\label{bmtherm}
BM_1=\frac{2\sigma_L^2(m+n-1)}{m^2n^2}\,,\,\,\,\,
BM_2=\frac{2\sigma_L^2\mu_D^2}{m^2\sigma_D^2}\,,\,\,\,\,
BM_3=\frac{2\mu_L^2}{n^2}+\frac{2\mu_L}{n}\,.
\end{equation}
In the tables below we have investigated the dependence of the bullwhip effect on the values of $m$ (the number of past lead times used in forecasting) for a given value of $n=5,10,20,30$ (the number of past demands used in forecasting), $\sigma_D/\mu_D=0.5$ (coefficient of demand variation), $\mu_L=3$ (expected value of lead times) and $\sigma_L=2$ (standard deviation of lead times). The tables show the impact of lead time forecasting on the bullwhip effect. It is evident that for small $m$
(e.g. $m=3$ or 5) the terms $BM_1$ and $BM_2$ contribute very much to the bullwhip effect and
when $m$ is large the impact of lead time forecasting on the bullwhip effect almost disappears but 
the effect remains by demand forecasting. For example for $n=5$ if $m$ changes from 3 to 50
the bullwhip effect measure varies from 6.72444 to 2.93971 (see Tab. \ref{bmm0}). This means that lead time forecasting can reduce the effect more than twice as much. Similarly e.g. in Tab. \ref{bmm3} if $m$ changes from 3 to 50
the bullwhip effect measure varies from 4.80716 to 1.23308 which indicates that the reduction in the effect can be almost four times as much. Moreover the effect is very big if $m=1$. This
follows from the fact that the forecast is based on the last known value of the lead time and the environment is ,,very random'' because we assume that lead times are mutually independent.

\begin{figure}[!H]
\begin{center}
\includegraphics[height=11cm,width=11cm]{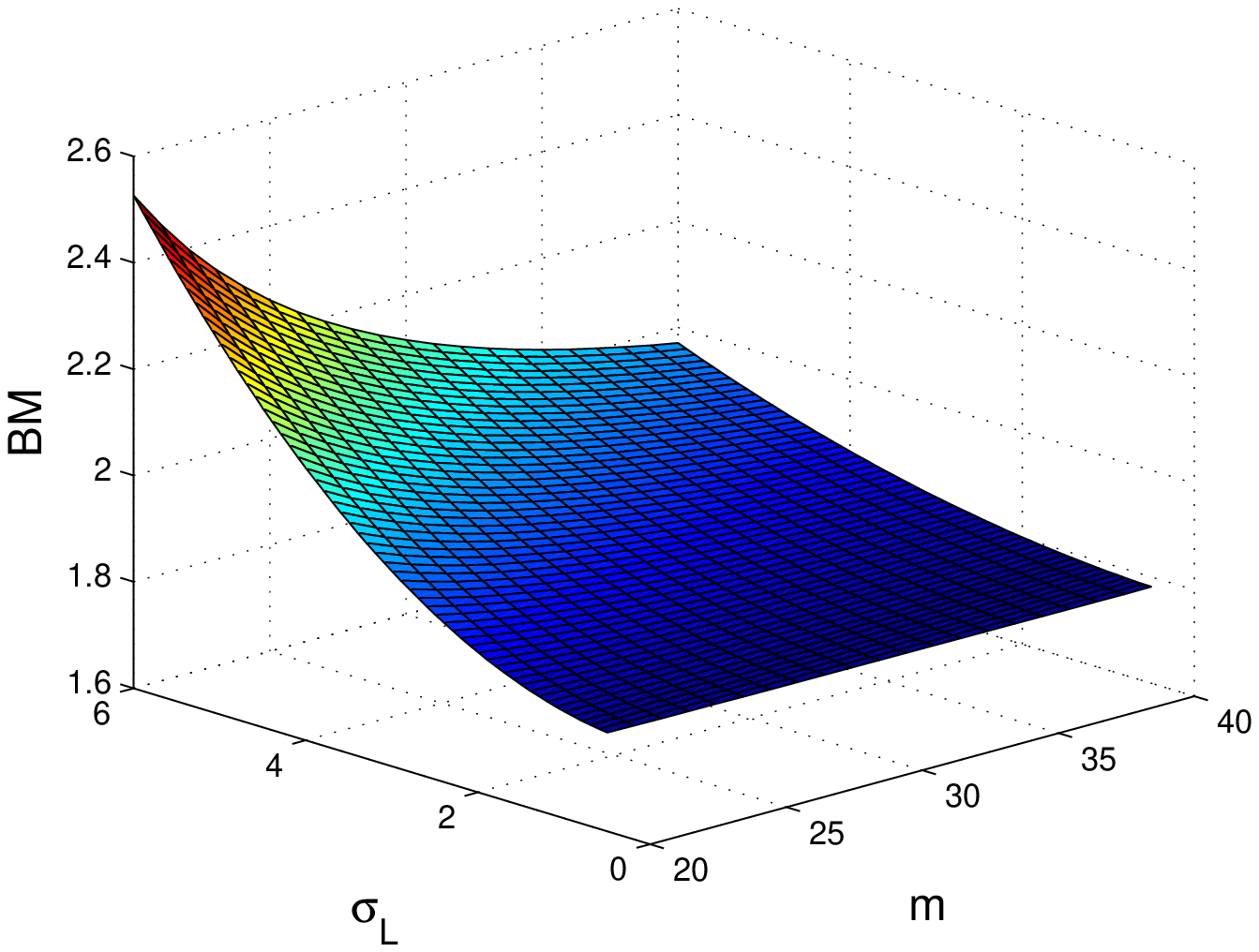}
\caption{The plot of the bullwhip effect measure as a function of $m$ and $\sigma_L$ where $n=10$, $\sigma_D/\mu_D=0.5$ and $\mu_L=3$.}\label{bmmsigl}
\end{center}
\end{figure}

\begin{figure}[!H]
\begin{center}
\includegraphics[height=11cm,width=11cm]{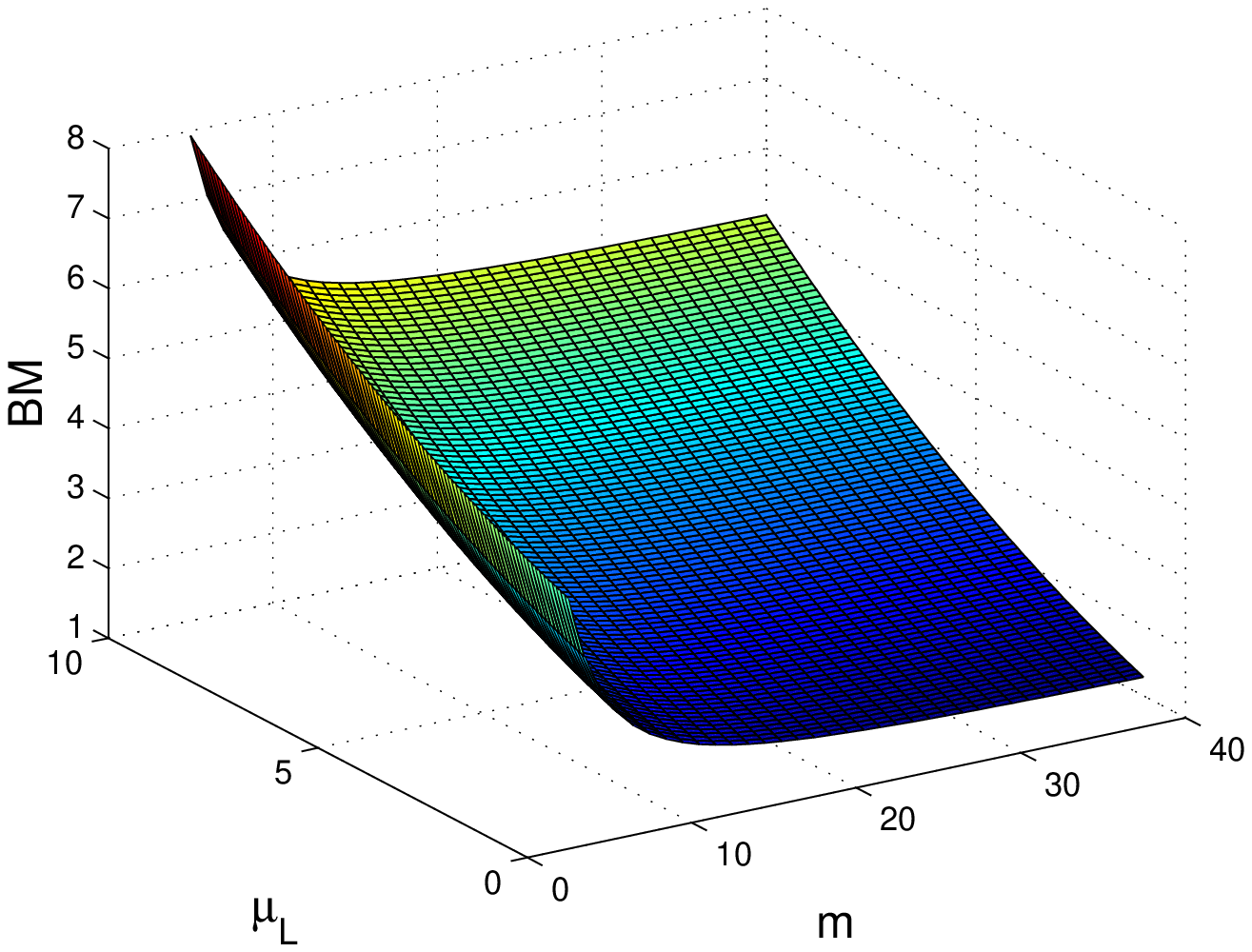}
\caption{The plot of the bullwhip effect measure as a function of $m$ and $\mu_L$ where $n=10$, $\sigma_D/\mu_D=0.5$ and $\sigma_L=3$.}\label{bmmmuL}
\end{center}
\end{figure}

\begin{figure}[!H]
\begin{center}
\includegraphics[height=11cm,width=11cm]{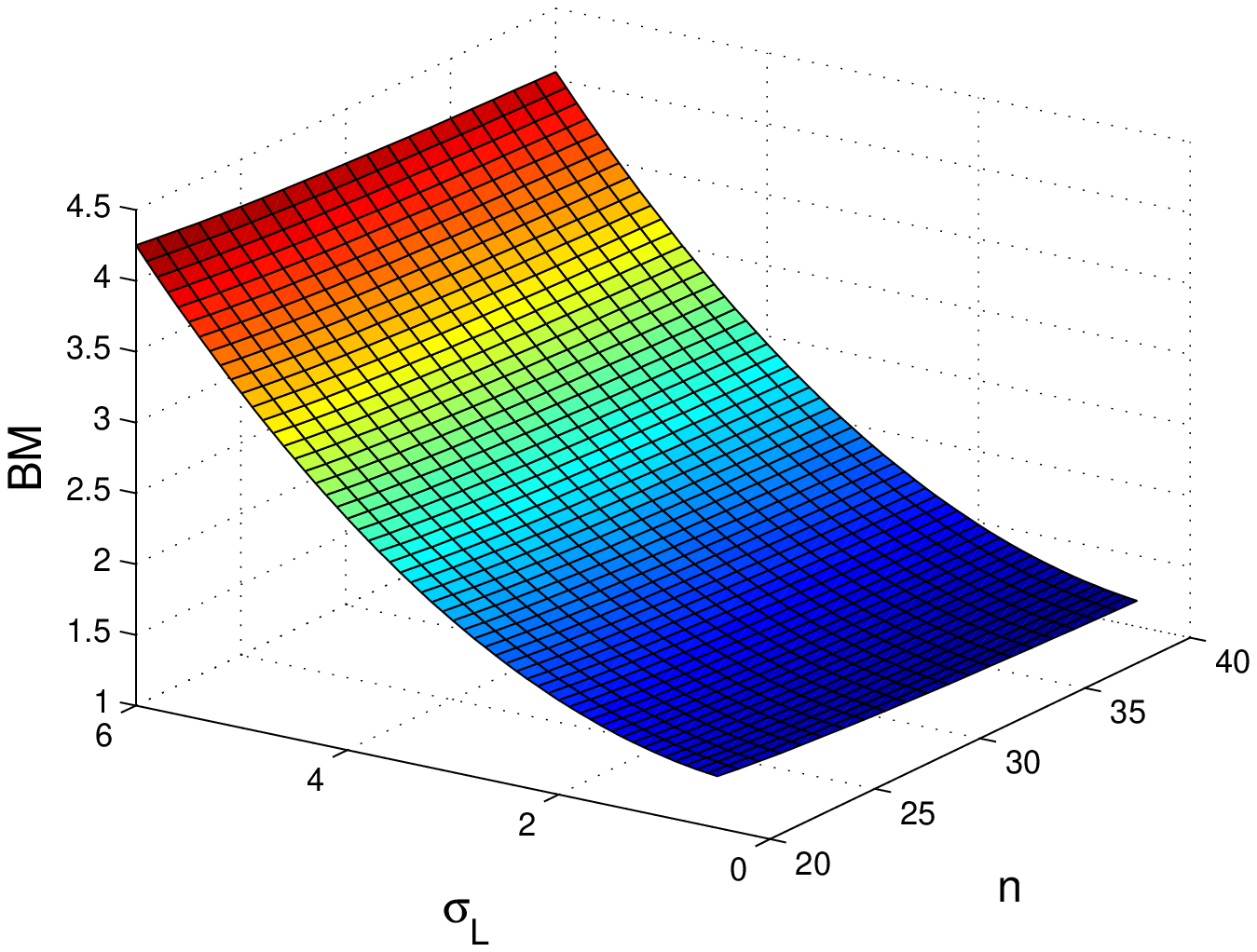}
\caption{The plot of the bullwhip effect measure as a function of $n$ and $\sigma_L$ where  $m=10$, $\sigma_D/\mu_D=0.5$ and $\mu_L=3$.}\label{bmnsigL}
\end{center}
\end{figure}

\begin{figure}[!H]
\begin{center}
\includegraphics[height=11cm,width=11cm]{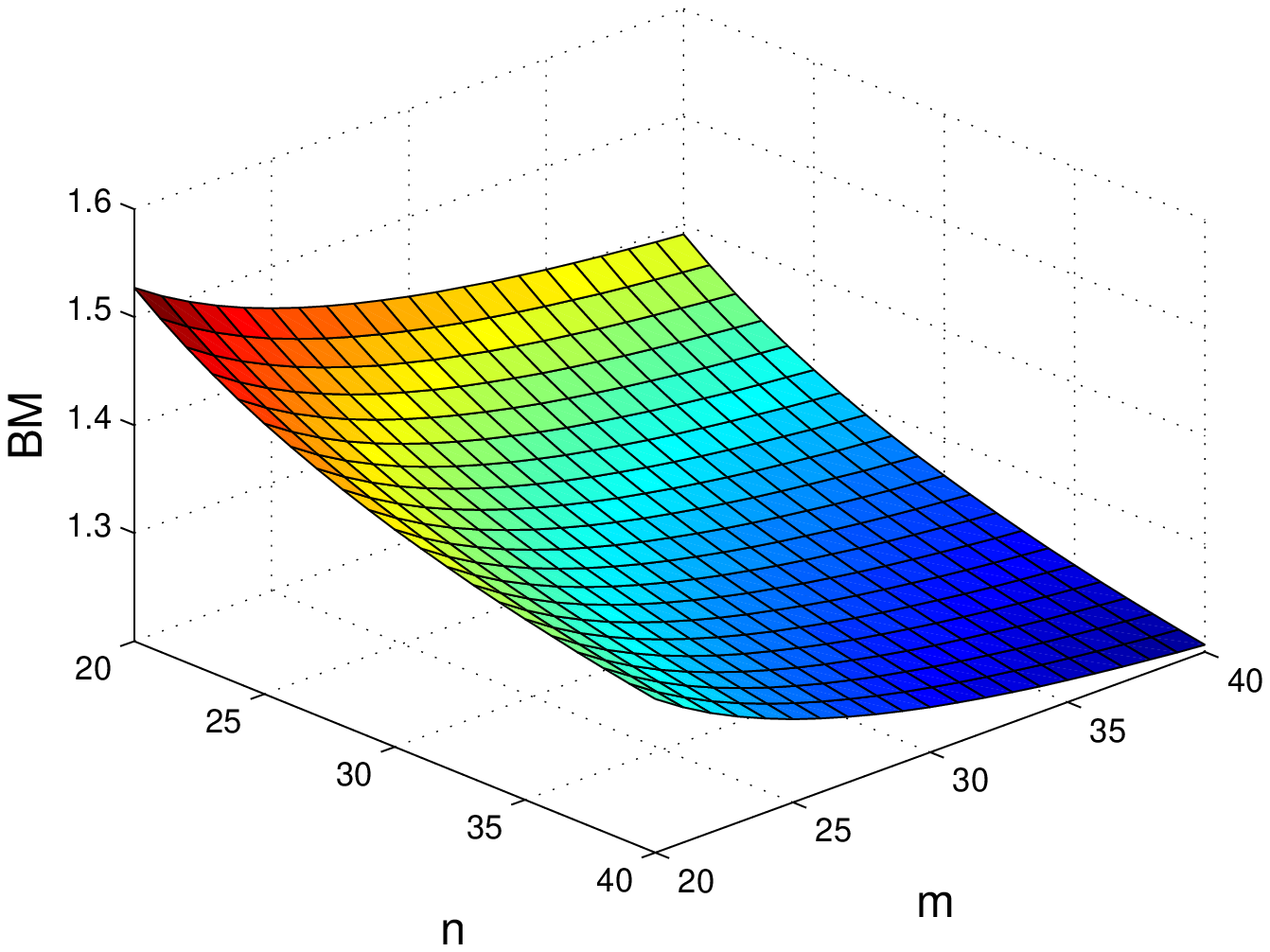}
\caption{The plot of the bullwhip effect measure as a function of $m$ and $n$ where $\sigma_D/\mu_D=0.5$, $\mu_L=3$ and $\sigma_L=3$.}\label{bmnm}
\end{center}
\end{figure}

\begin{figure}[!H]
\begin{center}
\includegraphics[height=11cm,width=11cm]{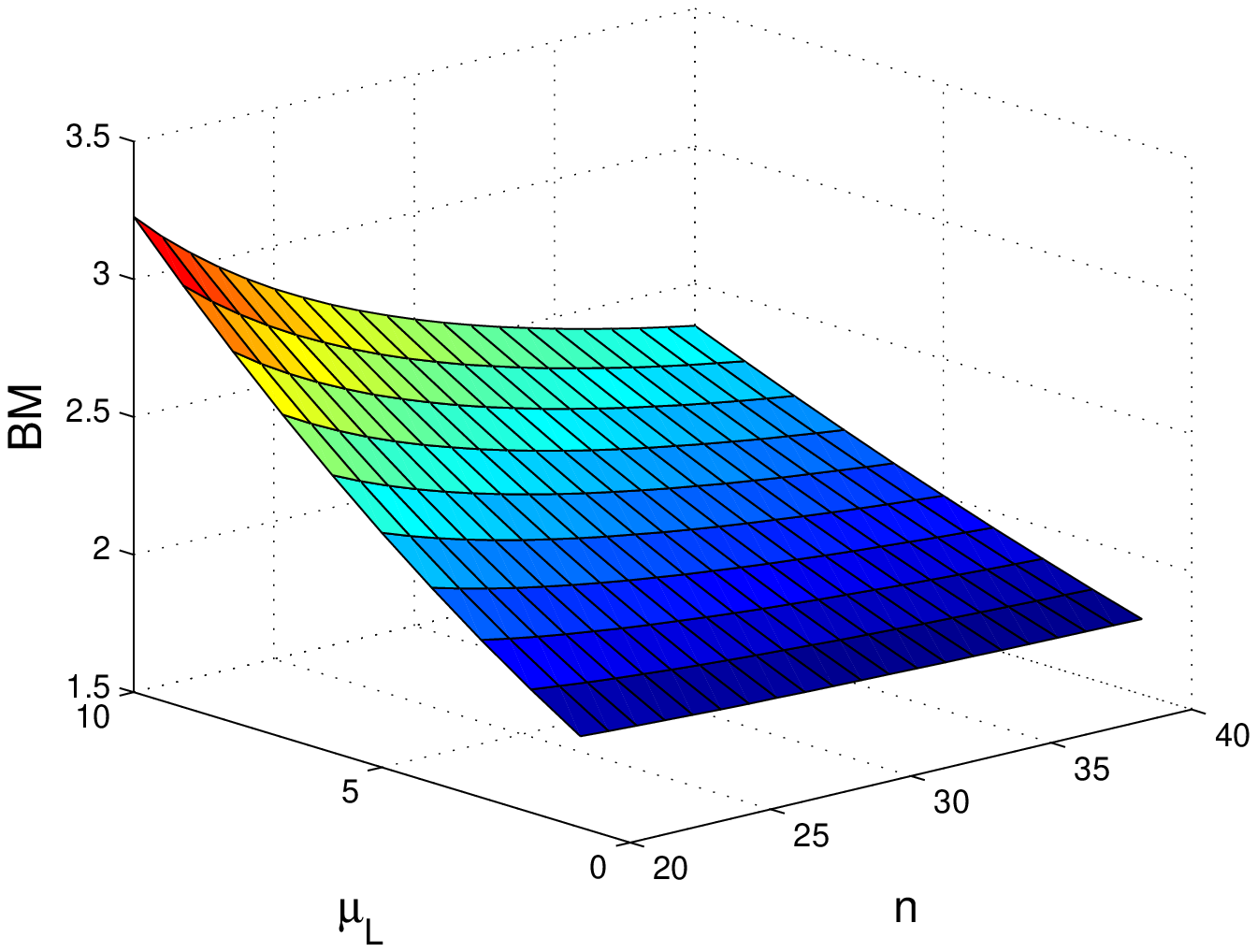}
\caption{The plot of the bullwhip effect measure as a function of $n$ and $\mu_L$ where $m=10$, $\sigma_D/\mu_D=0.5$ and $\sigma_L=3$.}\label{bmnmuL}
\end{center}
\end{figure}

We have also visualized the bullwhip effect measure as a function of two variables.
In Fig. \ref{bmmsigl} the measure of the bullwhip effect as a function of $m$ and $\sigma_L$ has been plotted where $m$ changes between 20 and 40, $\sigma_L\in [0.5, 6]$ and  $n=10$, $\sigma_D/\mu_D=0.5$ and $\mu_L=3$. 
The measure of the bullwhip effect as a function of $m=5,6,\ldots,40$ and $\mu_L\in[1, 10]$
for  $n=10$, $\sigma_D/\mu_D=0.5$ and $\sigma_L=3$ is shown in Fig. \ref{bmmmuL}. 
Fig. \ref{bmnsigL} shows the measure of the bullwhip effect depending on $n=20,21,\ldots,40$ and $\sigma_L\in [0.5, 6]$ for $m=10$, $\sigma_D/\mu_D=0.5$ and $\mu_L=3$. Similarly Fig. \ref{bmnm} presents the bullwhip effect as a function of $m$ and $n$ where their values change between 20 and 40 for $\sigma_D/\mu_D=0.5$, $\mu_L=3$ and $\sigma_L=3$. In Fig. \ref{bmnmuL} we visualize the measure of the bullwhip effect depending on $n=20,21,\ldots 40$ and $\mu_L\in [1, 10]$ where
$m=10$, $\sigma_D/\mu_D=0.5$ and $\sigma_L=3$.

\section {Conclusions and further research opportunities}
In this paper we have investigated the impact of lead time forecasting on the variance amplification in 
a simple two-stage supply chain with one supplier and one retailer, who employs the base stock policy for replenishment and the moving averages method for lead time and demand forecasting.
The exact form of the bullwhip effect measure indicates that lead time forecasting is a crucially contributing factor to the effect. The forecast of lead times gives two new summands $BM_1$ and $BM_2$ (see eq. (\ref{bmtherm})) in 
the bullwhip effect measure which substantially increase the value of the effect. These two summands are linear as a function
of lead time variance $\sigma_L^2$ and intensify the effect by the increase of $\sigma_L^2$. 
The next factor caused by lead time forecasting is the length of the sample of lead times used in the
forecast that is the value $m$ (see Th. \ref{bmmt}). If this value increases and goes to infinity then the variance amplification decreases and the impact of lead time forecasting disappears.
We should also note that the term $BM_1$ is of order $1/m$ for large $m$ that is $O(1/m)$
and the term $BM_2$ is of order $1/m^2$ for large $m$ that is $O(1/m^2)$ which means that the summand $BM_1$ has a bigger influence on the effect for large $m$. It is interesting that the term $BM_1$ can
be neglected if the length of demand observations applied in demand forecasting that is $n$ will be large because the summand $BM_1$ is also $O(1/n)$. Summarizing we ought to state that lead time forecasting is in fact a critically contributing factor to the bullwhip effect and its impact cannot be omitted in the design and management of supply chains. It is also worth noting that the effect stems from the need to estimate the lead time, and does not only depend on the expectation and variance of the lead time. 

The future research opportunities are widespread and necessary for the development of the supply chain management. In the further approaches to lead time forecasting problem we need to investigate 
other structures than iid of lead times and demands. Even if we consider more complicated structure of demands for example autoregressive-moving average leaving iid structure of lead times then this will complicate derivations of the bullwhip effect measure to a significant degree. Other opportunities lie in different forecasting methods for lead times and demands as well. Here we can apply different methods
for lead time forecasting and demand forecasting in a certain model or the same methods but other than the moving average method. In other directions of research one can investigate multi-echelon supply chains in the presence of stochastic lead times being predicted at every stage where the information on demands and lead times is shared or is not shared among the members of a supply chain. The value of the bullwhip effect measure in those situations 
will be surely valuable for theorists and  practitioners in the field of the supply chain management.

\begin{table}[!H]
  \begin{center}
  \caption{The measure of the bullwhip effect for $n=5$, $\sigma_D/\mu_D=0.5$, $\mu_L=3$ and $\sigma_L=2$ ($BM_3=1.920$).}\label{bmm0}
\vspace{2mm}
   \begin{tabular}{|c|c|c|c|}
     \hline
m & $BM_1$ & $BM_2$ & BM
\\ \hline
1   & 1.60000 & 32.00000 & 36.52000\\
3   & 0.24888 & 3.55555 & 6.72444\\
5   & 0.11520 & 1.28000 & 4.31520\\
10 & 0.04480 & 0.32000 & 3.28480\\
15 & 0.02702 & 0.14222 & 3.08924\\
20 & 0.01920 & 0.08000 & 3.01920\\
25 & 0.01484 & 0.05120 & 2.98604\\
30 & 0.01208 & 0.03555 & 2.96764\\
35 & 0.01018 & 0.02612 & 2.95631\\
40 & 0.00880 & 0.02000 & 2.94880\\
45 & 0.00774 & 0.01580 & 2.94354\\
50 & 0.00691 & 0.01280 & 2.93971\\
\hline
\end{tabular}
\end{center}
\end{table}

\begin{table}[!H]
  \begin{center}
  \caption{The measure of the bullwhip effect for $n=10$, $\sigma_D/\mu_D=0.5$, $\mu_L=3$ and $\sigma_L=2$ ($BM_3=0.780$).}\label{bmm1}
\vspace{2mm}
   \begin{tabular}{|c|c|c|c|}
     \hline
m & $BM_1$ & $BM_2$ & BM
\\ \hline
1   & 0.80000 & 32.00000 & 34.58000\\
3   & 0.10666 & 3.55555 & 5.44222\\
5   & 0.04480 & 1.28000 & 3.10480\\
10 & 0.01520 & 0.32000 & 2.11520\\
15 & 0.00853 & 0.14222 & 1.93075\\
20 & 0.00580 & 0.08000 & 1.86580\\
25 & 0.00435 & 0.05120 & 1.83555\\
30 & 0.00346 & 0.03555 & 1.81902\\
35 & 0.00287 & 0.02612 & 1.80899\\
40 & 0.00245 & 0.02000 & 1.80245\\
45 & 0.00213 & 0.01580 & 1.79793\\
50 & 0.00188 & 0.01280 & 1.79468\\
\hline
\end{tabular}
\end{center}
\end{table}

\begin{table}[!H]
  \begin{center}
  \caption{The measure of the bullwhip effect for $n=20$, $\sigma_D/\mu_D=0.5$, $\mu_L=3$ and $\sigma_L=2$ ($BM_3=0.345$).}\label{bmm2}
\vspace{2mm}
   \begin{tabular}{|c|c|c|c|}
     \hline
m & $BM_1$ & $BM_2$ & BM
\\ \hline
1   & 0.40000 & 32.00000 &33.74500 \\
3   & 0.04888 & 3.55555 & 4.94944\\
5   & 0.01920 & 1.28000 & 2.64420\\
10 & 0.00580 & 0.32000 & 1.67080\\
15 & 0.00302 & 0.14222 & 1.49024\\
20 & 0.00195 & 0.08000 & 1.42695\\
25 & 0.00140 & 0.05120 & 1.39760\\
30 & 0.00108 & 0.03555 & 1.38164\\
35 & 0.00088 & 0.02612 & 1.37200\\
40 & 0.00073 & 0.02000 & 1.36573\\
45 & 0.00063 & 0.01580 & 1.36143\\
50 & 0.00055 & 0.01280 & 1.35835\\
\hline
\end{tabular}
\end{center}
\end{table}

\begin{table}[!H]
  \begin{center}
  \caption{The measure of the bullwhip effect for $n=30$, $\sigma_D/\mu_D=0.5$, $\mu_L=3$ and $\sigma_L=2$ ($BM_3=0.220$).}\label{bmm3}
\vspace{2mm}
   \begin{tabular}{|c|c|c|c|}
     \hline
m & $BM_1$ & $BM_2$ & BM
\\ \hline
1   & 0.26666 & 32.00000 & 33.48666\\
3   & 0.03160 & 3.55555 & 4.80716\\
5   & 0.01208 & 1.28000 & 2.51208\\
10 & 0.00346 & 0.32000 & 1.54346\\
15 & 0.00173 & 0.14222 & 1.36396\\
20 & 0.00108 & 0.08000 & 1.30108\\
25 & 0.00076 & 0.05120 & 1.27196\\
30 & 0.00058 & 0.03555 & 1.25613\\
35 & 0.00046 & 0.02612 & 1.24658\\
40 & 0.00038 & 0.02000 & 1.24038\\
45 & 0.00032 & 0.01580 & 1.23612\\
50 & 0.00028 & 0.01280 & 1.23308\\
\hline
\end{tabular}
\end{center}
\end{table}

\subsection*{Acknowledgments} This work has been supported by the National Science Centre grant  \\  2012/07/B//HS4/00702.

\end{document}